\documentclass[a4paper,12pt]{article}
\usepackage{amsfonts,amssymb, pdfpages, xcolor}

\title{A numerical approximation for the standard one pressure system of two fluid flows with energy equations}
\author{M. Colombeau,\\ 
 Instituto de Matem\'atica, Estat\'istica e Computa\c c\~ao Cient\'ifica\\ Universidade Estadual de Campinas, Brazil.}

\begin{document}
\maketitle

\begin {abstract} We study numerically the standard one pressure model of two fluid flows with energy equations. This system is not solved  in time derivative. It has been transformed into an equivalent system solved in time derivative [S.T. Munkejord, S. Evje, T. Flatten, SIAM J. Sci. Comput. 31,4,2009,2587-2622]. We show that the scheme in this paper applies to both solved and nonsolved systems and gives same results. One usually adds a nonphysical term to render the system hyperbolic. However, explicit solutions and well posedness of the Cauchy problem for some nonlinear nonhyperbolic systems of physics have been obtained in some events by [B. Keyfitz et al.].  We also show that our scheme applies equally well to both  versions, with and without the additional term, whether solved in time derivative or not, which provides four versions of the system.  We observe that the nonhyperbolic and the hyperbolic systems give very close but slightly different results: the step values are always the same but peaks in gas and liquid velocities are observed in the nonhyperbolic model, which is typically observed in experimental results concerning the gas kick phenomenon though we are unable to say if this result is related or not. The numerical quality of the hyperbolic solved in time derivative system is better, therefore our results (same pressure and temperatures, and same step values in volume fraction and velocities besides the isolated peaks) provide a justification of the additional term that renders it hyperbolic.  It has also been observed that additional viscosity does not modify significantly the results for all four versions. Another difficulty lies in that these systems  are in nonconservative form and therefore its discontinuous solutions cannot make sense in the theory of distributions and it has been observed that different numerical schemes can lead to  different discontinuous solutions. For the hyperbolic system this solution  is identical to the main one in [S.T. Munkejord, S. Evje, T. Flatten] obtained from completely different methods. Finally, we improve the order one scheme in use above by transforming it explicitly into  order three in space.\\

\end{abstract}
e-mail: m.colombeau@orange.fr\\
\textit{*this research has been  done thanks to  financial support of FAPESP, processo 2012/15780-9, when the author was doing her postdoctoral studies at the University of S\~ao Paulo, USP, S\~ao Paulo, Brazil.}\\


\textbf{1.  Introduction}.\\ 
We study numerically the  basic system  used to model a mixture of two immiscible fluids from the conservation laws of mass, momentum and energy with the natural assumption that the pressures at the same point are equal inside the two fluids \cite{Wendroff} p. 373, \cite{Munk} p. 2589 and 2590, \cite{FFMunk} p. 413, \cite{FlattenRoe} p. 479, \cite{EvjeFlatten}  p. 179, \cite{Cortes} p. 465,\dots. This paper contributes to the research on this system in five ways.\\

$\bullet$ First, since the system is in nonconservative form one should state precisely inside the schemes the formulas of discretization of the nonconservative products since various slight changes could lead to different weak solutions \cite{Munk} p. 2620. We observe that in the case we admit that the jumps between successive cells are "`small"' (which permits large jumps in the numerical solution provided a number of cells are involved inside these jumps) the ambiguity can be solved in a natural way by assuming that the step functions (constant on each cell) at times $n\Delta t$  satisfy the equations. One recovers the simplest formulas proposed by various authors but we have a mathematical justification of these formulas. \\

$\bullet$ Second, a sequence of approximate solutions with full mathematical proof that they tend to satisfy the system has been constructed in \cite{Colombeautheo} to play the role of a (lacking) explicit exact solution. It has been checked in \cite{Colombeautheo} that the scheme in this paper produces exactly this theoretical solution in both cases ($\delta=0$ and $\delta=2$, see below). This shows that the scheme in sections 3 and 4 gives a solution of the system but, in absence of a theoretical uniqueness result, one cannot claim this is the physically correct solution, although there is an accumulation of indications at least for the Toumi shock tube problem.\\

$\bullet$ Our scheme permits to treat as well the original system  \cite{Wendroff} p. 373, \cite{Munk} pp.2589-2590 and p. 2592, not solved in time derivatives and the system solved in time derivatives obtained in \cite{Munk} p. 2595-2596. We check that the two systems give exactly the same numerical solution, which is not trivial in the case of shock waves since nonlinear calculations that preserve the smooth solutions do not in general preserve the shock waves: it is well known that the equations $u_t+uu_x=0$ and $uu_t+u^2u_x=0$ have different jump conditions.\\

$\bullet$ The authors on multifluid flows consider an additional term to render the system hyperbolic, then study the hyperbolic system so obtained. Motivated by the fact that B. Keyfitz et al \cite{Keyfitzlivro,Keyfitzart1,Keyfitzart2,Keyfitzart3,KeyfitzLopes,KeyfitzSan,Keyfitz0,Keyfitzrev1,Keyfitzrev2,Keyfitzrev3} have discovered that, in certain cases, nonlinear nonhyperbolic systems could have realistic physical solutions we study also numerically the original nonhyperbolic system. We observe same step values but a significant difference in form of a well defined peak in liquid flow velocity and gas flow velocity present in the original system and absent after the additional term. Comparison with some experimental data suggest that these peaks could be the gas kick phenomenon.\\

$\bullet$ Finally we transform our scheme into higher order schemes in view of an extension to multidimension.\\

Now we recall the equations. The two fluids are denoted by the indices $1$  and $2$, for instance mixture of oil and natural gas in extraction tubes of oil exploitation for the gas kick simulation \cite{Avelar,Bendiksen,Larsen}. The system, as it stems from physics, called canonical nonconservative form, is stated as\\

\begin{equation}\frac{\partial}{\partial t}(\rho_1 \alpha_1)+\frac{\partial}{\partial x}(\rho_1\alpha_1 v_1)=0,\end{equation}
\begin{equation}\frac{\partial}{\partial t}(\rho_2 \alpha_2)+\frac{\partial}{\partial x}(\rho_2\alpha_2 v_2)=0,\end{equation}
\begin{equation}\frac{\partial}{\partial t}(\rho_1 \alpha_1v_1)+\frac{\partial}{\partial x}(\rho_1\alpha_1( v_1)^2)+\alpha_1\frac{\partial p}{\partial x}+\tau_i=g\alpha_1\rho_1,\end{equation}
\begin{equation}\frac{\partial}{\partial t}(\rho_2 \alpha_2v_2)+\frac{\partial}{\partial x}(\rho_2\alpha_2( v_2)^2)+\alpha_2\frac{\partial p}{\partial x}-\tau_i=g\alpha_2\rho_2,\end{equation}

\begin{equation}\frac{\partial}{\partial t}(E_1)+p\frac{\partial}{\partial t}(\alpha_1) +\frac{\partial}{\partial x}(E_1v_1+p\alpha_1v_1)+v_{\tau}\tau_i=g\rho_1\alpha_1v_1,   \end{equation}
\begin{equation}\frac{\partial}{\partial t}(E_2)+p\frac{\partial}{\partial t}(\alpha_2) +\frac{\partial}{\partial x}(E_2v_2+p\alpha_2v_2)-v_{\tau}\tau_i=g\rho_2\alpha_2v_2,   \end{equation}
\\
which are the conservation laws of mass, momentum and energy for fluids 1,2. They are complemented by

\begin{equation}\alpha_1+\alpha_2=1,\end{equation}
\begin{equation} \tau_i=\Delta p\frac{\partial\alpha_1}{\partial x}, \ \ \Delta p=\delta\frac{\alpha_1\alpha_2\rho_1\rho_2}{\rho_1\alpha_2+\rho_2\alpha_1}(v_1-v_2)^2,\end{equation}
\begin{equation} v_\tau=\frac{\alpha_2\gamma_1v_1+\alpha_1\gamma_2v_2}{\alpha_2\gamma_1+\alpha_1\gamma_2},\end{equation}
\\
$\delta\geq 0, \gamma_i\geq 0$ and the state laws, 
\begin{equation}p_1=(K_1-1)(\frac{E_1}{\alpha_1}-\frac{\rho_1 (v_1)^2}{2})-K_1p_{\infty,1},\end{equation} \begin{equation}p_2=(K_2-1)(\frac{E_2}{\alpha_2}-\frac{\rho_2 (v_2)^2}{2})-K_2p_{\infty,2},\end{equation}
in which we postulate $p_1=p_2$, denoted by $p$, from the equal pressure assumption.\\

 The physical variables are 
the densities $\rho_i(x,t)$, the velocities $v_i(x,t)$, the volumic proportions
$\alpha_i(x,t)$, the equal pressures
$p_i(x,t)$,   the total energy densities $E_i(x,t)$, $i=1,2$. The constant
$g$ is the component of the gravitational acceleration in the direction of the tube,
$\tau_i$ is the interphasic momentum exchange term, $\Delta p$ is an interface pressure correction term
and $v_\tau$ is an operator with the dimension of a velocity which stems from the laws of thermodynamics
as exposed in \cite{Munk} where the choice of the formulas (8,9) is justified. The letters $\delta, \gamma_1, \gamma_2, K_1, K_2, p_{\infty,1}$ and $p_{\infty,2}$ represent real values which are chosen and explained in \cite{Munk,FFMunk}. We will adopt in the numerical calculations the values used in \cite{Munk} p.2613 and \cite{FFMunk} p.346 for the Toumi shock tube problem. \\ 

The system presents three basic peculiarities: obviously it is not solved in the time derivative and it is not in conservation form; non obviously  it is not hyperbolic \cite{Munk} p. 2595. Each of these three peculiarities already causes serious problems for the elaboration of numerical schemes. In \cite{Munk} pp. 2592-2596 these authors, by a deep analysis  involving the laws of thermodynamics, have produced a formally equivalent formulation which is resolved in the time derivative and has served to construct numerical methods, \cite{Munk,FFMunk}. The nonconservative character implies a  nonuniqueness of stable (entropic) solutions \cite{Munk} pp. 2587,2620, with the problem of searching the correct ones. The absence of hyperbolicity generally suggests absence of well posedness. On the other hand this system is of basic importance in industry such as oil extraction in deep sea \cite{Avelar,Bendiksen,Larsen} or cooling of nuclear power stations \cite{Bestion,WAHA3}. In this paper we will address these three peculiarities.\\

Concerning nonhyperbolicity, somewhat unexpectedly,   B. Keyfitz has pointed out various instances in which nonlinear nonhyperbolic systems modeling physics  have  physically meaningful solutions with applications in various domains, for instance in traffic flow, multifluid flows and porous medium, besides the fact the corresponding  linear systems could be ill posed, see \cite{Keyfitzlivro} pp. 150-151, the articles \cite{Keyfitzart1,Keyfitzart2,Keyfitzart3,KeyfitzLopes,KeyfitzSan} and the review papers \cite{Keyfitz0,Keyfitzrev1,Keyfitzrev2,Keyfitzrev3}. This could be explained by the fact that some positive nonlinear effects could be lost in the linearization. This motivates an attempt  for a numerical study of the nonhyperbolic system, which has been shown realistic in the four equations model \cite{ColombeauJDE}.\\

The scheme that we propose has been constructed from a scheme proposed in \cite{ColombeauSiam,ColombeauNMPDE} for pressureless fluids, \cite{ColombeauJCAM} for ideal gases and shallow water, and \cite{ColombeauJCAM2} for various systems. A version for the four equations model of two fluid flows is given in \cite{ColombeauJDE} section 7. We begin to expose a scheme which is of order one in time and space. Then we extend it to a scheme of order three in space that gives more accurate results.     The scheme  applies 
 to system (1-11), called  canonical nonconservative form,  although it is not solved in the time derivative, not conservative and not hyperbolic. Applying it with the  additional term (8) used in \cite{Munk} ($\delta=2$) to render the system hyperbolic one observes (figure 1) that one obtains exactly the numerical results in \cite{Munk} p. 2615,2617 and \cite{FFMunk}  p. 437.\\

Using  the  systems solved and nonsolved in the time derivatives \cite{Munk}  we obtain  (figures 1 and 2) with the scheme in this paper the same result for the discontinuous solution of the Toumi problem, which constitutes  a numerical verification of the scheme and a numerical verification of the transformation of the equations in the case of discontinuous solutions \cite{Munk} pp. 2592-2596.  Another interesting point is that the possibility to use the scheme with and without the additional term ($\delta=2$ and $\delta=0$ respectively in (8), figure 3) will permit to observe the influence of additional terms on the numerical solution and also to study numerically the original nonhyperbolic model \cite{Wendroff} p.373. These observations will be discussed in the paper and will finally provide a justification of the relevance of the additional term.\\

A theorem in section 4 shows, both for the hyperbolic and the nonhyperbolic model, that the approximate 
solutions from the scheme tend to satisfy the equations somewhat independently of a  reasonable arbitrariness in a definition of a nonconservative product needed to give a sense to the equations.\\

\textbf{2. A choice of the nonconservative products inside the formulas of the scheme.} A basic problem lies in the nonconservative form of the formulas of the equations for example the single term $\alpha_i\partial_xp$ in (3,4). It is known that different numerical schemes can converge to different weak solutions \cite{Munk} p. 2620. Therefore it is of basic importance to decide a precise choice of the formulas inside a numerical scheme for the discretization of the nonconservative products. In this section we propose such a choice under a "`small jump assumption"' commonly used in numerical schemes: we postulate that the jumps $s_i$ between a cell $[(i-\frac{1}{2})h, (i+\frac{1}{2})h]$  and its neighbour cells are such that the higher order powers $(s_i)^n, \ n=2,3$ can be neglected in the calculations,  in which we retain therefore only the first order terms $s_i$. \\

Then the  reasoning that leads to a precise formula inside the scheme for the nonconservative products is as follows.  At times $t_n=n\Delta t$ one would like to have a numerical solution in form of step functions on the cells that would satisfy the equations as much as possible. How are these step functions in the case of a nonconservative product under the small jumps assumption?  Consider a general nonlinear equation 
\begin{equation} \partial_t u+A(u,v)\partial_x u+B(u,v)\partial_x v=0\end{equation}
where $A$ and $B$ are polynomials in $u,v$. With standard notation, at the interface $(i+\frac{1}{2})h$, we set $u,v$ of the form
\begin{equation} u(x,t,h)=u_i+(u_{i+1}-u_i)H_u(x-ct), \ v(x,t,h)=v_i+(v_{i+1}-v_i)H_v(x-ct), \end{equation}
where $H_u$ and $H_v$ are regularizations of the Heaviside functions $H$: here $H$ is discontinuous at $(i+\frac{1}{2})h$ and $H_u, H_v$ are smooth functions that  jump from 0 to 1 in the interval $]ih,(i+1)h[$. Then with this mollification the nonconservative products make sense. What is their value in the mollified case, to be adopted at the limit of absence of mollification? To this end plug the formulas (13) into the equation (12). One obtains 
\begin{equation} -c(u_{i+1}-u_i)H_u'+A(u,v)(u_{i+1}-u_i)H_u'+B(u,v)(v_{i+1}-v_i)H_v'=0.\end{equation}
From the small jump assumption one replaces $A(u,v) $ and $B(u,v)$ by $A(u_i,v_i) $ and $B(u_i,v_i)$ respectively. From (14) by integration one deduces that $H_u=H_v$ since both jump from 0 to 1 and are proportional. Denoting $\overline{H}=H_u=H_v$ a product such as $uv_x$ gives $$\int_ih^{(i+1)h}uv_x dx=\int [u_i+(u_{i+1}-u_i)\overline{H}](v_{i+1}-v_i)\overline{H}'dx=u_i(v_{i+1}-v_i)+(u_{i+1}-u_i)(v_{i+1}-v_i)\frac{1}{2}$$ $=\frac{u_{i+1}-u_i}{2}(v_{i+1}-v_i).$\\
\\
Therefore one is led, under the small jump assumption, to discretize $uv_x$ by \begin{equation}\frac{1}{h}\frac{u_{i+1}+u_i}{2}(v_{i+1}-v_i)\end{equation}, or equivalently, seeking the terms in factor of $u_i$ for convenience, by \begin{equation}u_i( \frac{v_{i+1}-v_{i-1}}{2h})\end{equation} which is simply a familiar centered discretization (formulas (15) and (16) are equivalent since we will understand the right-hand side of the equations in the sense of distributions). This discretization will be applied to all nonconservative products in the scheme below. The unique point would be to check that there will appear only small values of jumps between successive cells even inside physical discontinuities such as shock waves.\\

According to \cite{Munk} p. 2597-2598 the above formula for nonconservative products consists in stating that the paths are same for $u$ and $v$ in the small jumps between two successive cells. Indeed the mollifications $H_u$ and $H_v$ create a path for $H_u$ and a path for $H_v$, not necessarily equal a priori. Then we prove under the small jumps assumption that these paths should be equal in order that the numerical solutions at times $t_n=n\Delta t$ stick to the equations. Since the jumps in $u$ and in $v$ between the various cells do not necessarily remain proportional throughout a shock wave, the global paths in $u$ and $v$ for a physical jump involving a number of cells inside the jump can be very different, even when the paths between successive cells are same (an axample of this fact is given in \cite{Abgrall}). \\

Note that due to their simplicity formulas (15) and (16) are widely used by various authors to discretize nonconservative products: for instance formula (16) is used in \cite{Abgrall} p.2761 between formulas (7) and (8) there. The novelty is that we provide a justification: under the small jump assumption other choices would give step functions (constant on each cell) at times $n\Delta t$ that would not be solutions of the equations at times $t_n$. The authors of \cite{Abgrall} give an example of a nonconservative system for which application of this formula for the jumps between successive cells  does not produce the same jump formulas for the global jumps in the shock waves obtained from the scheme. The scheme in this paper can also produce this natural phenomenon since the jumps between successive are in general not proportional for all interfaces of cells inside a physical jump.\\  


\textbf{ 3. A method of sequences of smooth approximate solutions with full mathematical proof.} The system (1-11) cannot have discontinuous solutions in the sense of distribution theory in case of shock waves because of the terms $\alpha_i \partial_x p$ and $p\partial_t\alpha_i$ which do not permit a transfer of the derivatives to smooth test functions. Therefore one has no other choice than approximation of the values of $\alpha_i$ and $p$ by smooth functions. Indeed this method has been widely developed in mathematics for systems without solutions in the sense of distributions and it has given  very neat results for conservative systems, such as  studies of $\delta$ waves and  $\delta^{(n)}$ waves, in \cite {Danilov1,DanilovM1,DanilovM2,DanilovOS,DanilovO1,DanilovO2,DanilovS1,DanilovS2,Panov,Shelkovich2,ShelkovichRMS} among other papers. In this section we recall results in \cite{Colombeautheo}.   For the nonconservative system (1-11) a sequence of smooth approximate solutions with a full proof that the sequence tends to satisfy the equations has been constructed in \cite{Colombeautheo}. This construction is based on a system of 6 ODEs in Banach space relevant of the Lipschitz theory of ODEs, therefore they can be solved numerically by convergent classical schemes for ODEs. The  defect of this approach lies in the present absence of a uniqueness result of a    "`limit"  (that could be simply a limit from observation) that "`admissible sequences"' could produce. The sequence constructed in \cite{Colombeautheo} has given exactly the results obtained in figure 2 from the schemes in sections 4. Of course one could argue that the construction of this sequence of smooth solutions has been inspired by this scheme-although noticeably different, with intervention of physical arguments such as the convolution in pressure.\\

\textbf{4. A numerical  scheme for the canonical nonconservative form (1-11).} The initial scheme in     \cite{ColombeauSiam} has been extended to the more complicated system of collisional self gravitating fluids \cite{ColombeauNMPDE, ColombeauJMP}, to the systems of ideal gases  and shallow water  equations \cite{ColombeauJCAM,ColombeauZeit}, to the singular shocks and the $\delta^{(n)}$ shocks of the sytems of Keyfitz-Kranzer \cite{Keyfitz1,Keyfitz2} and Panov-Shelkovich \cite{Panov}  in \cite{ColombeauJCAM2} and to the four equations model of two fluid flows in \cite{ColombeauJDE} section 7.  
The basic method  of the scheme    consists in a splitting of the equations (1-6) into a transport step: transport of the fluid i with velocity $v_i$, i=1,2 (though (16,17) are not exactly a transport). 

\begin{equation}\frac{\partial}{\partial t}(\rho_1 \alpha_1)+\frac{\partial}{\partial x}(\rho_1\alpha_1 v_1)=0,\end{equation}
\begin{equation}\frac{\partial}{\partial t}(\rho_2 \alpha_2)+\frac{\partial}{\partial x}(\rho_2\alpha_2 v_2)=0,\end{equation}
\begin{equation}\frac{\partial}{\partial t}(\rho_1 \alpha_1v_1)+\frac{\partial}{\partial x}(\rho_1\alpha_1( v_1)^2)=0,\end{equation}
\begin{equation}\frac{\partial}{\partial t}(\rho_2 \alpha_2v_2)+\frac{\partial}{\partial x}(\rho_2\alpha_2( v_2)^2)=0,\end{equation}

\begin{equation}\frac{\partial}{\partial t}(E_1)+p\frac{\partial}{\partial t}(\alpha_1) +\frac{\partial}{\partial x}((E_1+p\alpha_1)v_1)=0,   \end{equation}
\begin{equation}\frac{\partial}{\partial t}(E_2)+p\frac{\partial}{\partial t}(\alpha_2) +\frac{\partial}{\partial x}((E_2+p\alpha_2)v_2)=0,   \end{equation}
\\
and a pressure correction step

\begin{equation}\frac{\partial}{\partial t}(\rho_1 \alpha_1)=0,\end{equation}
\begin{equation}\frac{\partial}{\partial t}(\rho_2 \alpha_2)=0,\end{equation}
\begin{equation}\frac{\partial}{\partial t}(\rho_1 \alpha_1v_1)+\alpha_1\frac{\partial p}{\partial x}+\tau_i=g\alpha_1\rho_1,\end{equation}
\begin{equation}\frac{\partial}{\partial t}(\rho_2 \alpha_2v_2)+\alpha_2\frac{\partial p}{\partial x}-\tau_i=g\alpha_2\rho_2,\end{equation}

\begin{equation}\frac{\partial}{\partial t}(E_1)+p\frac{\partial}{\partial t}(\alpha_1) +v_{\tau}\tau_i=g\rho_1\alpha_1v_1,   \end{equation}
\begin{equation}\frac{\partial}{\partial t}(E_2)+p\frac{\partial}{\partial t}(\alpha_2) -v_{\tau}\tau_i=g\rho_2\alpha_2v_2,   \end{equation}

 \noindent  with, in between, an averaging step and, in the case of system (1-11), auxiliary  algebraic calculations not needed in the case of a single fluid. We notice the groups $\frac{\partial}{\partial t}(E_j)+p\frac{\partial}{\partial t}(\alpha_j), \ j=1,2$, which are not  time derivatives. We simply discretize them as 
\begin{equation}\frac{(E_j)(t+\Delta t)-(E_j)(t)+p(t)  ((\alpha_j)(t+\Delta t)-(\alpha_j)(t))}{\Delta t}.\end{equation}

 Now we give the formulas of the scheme.
 We set $r_j=\rho_j\alpha_j$ and $F_j=E_j+p\alpha_j, \ j=1,2$. The notation $F$ serves to contain $E$ and $\alpha$ as a whole for  time evolution according to formula (29), them to compute $E_j(t+\Delta t)$ and $\alpha_j(t+\Delta t)$ from algebraic calculations. We assume we have computed at time $t_n=n\Delta t$ the family
 \begin{equation}\{(r_1)_i^n,(r_2)_i^n,(r_1v_1)_i^n, (r_2v_2)_i^n,  (F_1)_i^n, (F_2)_i^n,p_i^n\}_{i\in\mathbb{Z}}.\end{equation} 
We seek the values at $t_{n+1}=(n+1)\Delta t$. We set as usual $r=\frac{\Delta t}{h}$ where $h$ is the space step. We set  $(v_j)_i^n=\frac{(r_jv_j)_i^n}{(r_j)_i^n}$, and $ ((v_j)_i^n)^+=max((v_j)_i^n,0), ((v_j)_i^n)^-=max(-(v_j)_i^n,0), j=1,2.$ \\

$\bullet$ \textbf{first step: transport.} 
\begin{equation}(\widetilde{r_1})_i=(r_1)_i^n+r[(r_1)_{i-1}^n((v_1)_{i-1}^n)^+-(r_1)_{i}^n(|v_1|_{i}^n)+(r_1)_{i+1}^n((v_1)_{i+1}^n)^-],\end{equation}
\begin{equation}(\widetilde{r_2})_i=(r_2)_i^n+r[(r_2)_{i-1}^n((v_2)_{i-1}^n)^+-(r_2)_{i}^n(|v_2|_{i}^n)+(r_2)_{i+1}^n((v_2)_{i+1}^n)^-],\end{equation}
\begin{equation}(\widetilde{r_1v_1})_i=(r_1v_1)_i^n+r[(r_1v_1)_{i-1}^n((v_1)_{i-1}^n)^+-(r_1v_1)_{i}^n(|v_1|_{i}^n)+(r_1v_1)_{i+1}^n((v_1)_{i+1}^n)^-],\end{equation}
\begin{equation}(\widetilde{r_2v_2})_i=(r_2v_2)_i^n+r[(r_2v_2)_{i-1}^n((v_2)_{i-1}^n)^+-(r_2v_2)_{i}^n(|v_2|_{i}^n)+(r_2v_2)_{i+1}^n((v_2)_{i+1}^n)^-],\end{equation}
\begin{equation}(\widetilde{F_1})_i=(F_1)_i^n+r[(F_1)_{i-1}^n((v_1)_{i-1}^n)^+-(F_1)_{i}^n(|v_1|_{i}^n)+(F_1)_{i+1}^n((v_1)_{i+1}^n)^-],\end{equation}
\begin{equation}(\widetilde{F_2})_i=(F_2)_i^n+r[(F_2)_{i-1}^n((v_2)_{i-1}^n)^+-(F_2)_{i}^n(|v_2|_{i}^n)+(F_2)_{i+1}^n((v_2)_{i+1}^n)^-],\end{equation}

For $\omega=r_1,r_2,r_1v_1$ and $r_2v_2$ \ $\widetilde{\omega}$ stands for $\omega(t+\Delta t)$ while $\widetilde{F_j}, j=1,2$ will be treated as $E_j(t+\Delta t) +p(t)\alpha_j(t+\Delta t)$ from which one will extract $E_j(t+\Delta t)$ and $\alpha_j (t+\Delta t)$ by algebraic calculations.\\

 Now we do  an averaging of the quantities calculated above. This averaging is indispensable in most cases as exposed in \cite{ColombeauNMPDE,ColombeauJCAM,ColombeauJCAM2}.\\

$\bullet$ \textbf{second step: averaging of the quantities calculated above.} Let $a, 0<a<\frac{1}{2}$ be a given real number. We replace each of the values generically denoted $(\widetilde{q})_i$ calculated in the first step (31-36) by $a(\widetilde{q})_{i-1}+(1-2  a)(\widetilde{q})_i+a(\widetilde{q})_{i+1}$. This last quantity is denoted  by $\overline{q}_i$. The choice of the value $a, \ 0<a<\frac{1}{2}$, is done from numerical tests as explained in \cite{ColombeauJCAM} p.17. Therefore after this averaging we have the six families of values
 \begin{equation}\{(\overline{r_1})_i, (\overline{r_2})_i, (\overline{r_1v_1})_i, (\overline{r_2v_2})_i,
(\overline{F_1})_i,(\overline{F_2})_i\}_{i\in\mathbb{Z}}.\end{equation}
Further, we set
\begin{equation}(\overline{v_1})_i=\frac{(\overline{r_1v_1})_i}{ (\overline{r_1})_i}, (\overline{v_2})_i=\frac{(\overline{r_2v_2})_i}{ (\overline{r_2})_i}.\end{equation}

$\bullet$ \textbf{third step: auxiliary algebraic calculations.} We compute values denoted by  $\overline{q}   $  for $q=E_1, E_2, \alpha_1, \alpha_2, p_1$ and $p_2$ from the values (37,38). This could also be done from the $\widetilde{q}$ values in (31-36), defining $v_1$ and $v_2$ by the usual quotient. We exploit the four algebraic formulas (10,11,7,$p_1=p_2$) and the two values 
$[(\overline{E_1})_i+p_i^n(\overline{\alpha_1})_i]=(\overline{F_1})_i$ and $[(\overline{E_2})_i+p_i^n(\overline{\alpha_2})_i]=(\overline{F_2})_i$  obtained above. This gives 6 equations for the 6 unknown values  $(\overline{E_1})_i, (\overline{E_2})_i, (\overline{\alpha_1})_i, (\overline{\alpha_2})_i, (\overline{p_1})_i $ and $(\overline{p_2})_i$. After some easy algebraic calculations one obtains that $(\overline{\alpha_1})_i$ is the root located in $]0,1[$ of the equation $AX^2+BX+C=0$ where
\begin{equation}A=(K_1-1)p_i^n+K_1p_{\infty,1}-(K_2-1)p_i^n-K_2p_{\infty,2},\end{equation}
$ B=-(K_1-1)(\overline{F_1})_i-(K_2-1)(\overline{F_2})_i-K_1p_{\infty,1}+K_2p_{\infty,2}+\frac{1}{2}(K_1-1)(\overline{r_1v_1})_i(\overline{v_1})_i+$\begin{equation}\frac{1}{2}(K_2-1)(\overline{r_2v_2})_i(\overline{v_2})_i+(K_2-K_1)p_i^n,\end{equation}
\begin{equation}C=(K_1-1)(\overline{F_1})_i-\frac{1}{2}(K_1-1)(\overline{r_1v_1})_i(\overline{v_1})_i.\end{equation}
Once one has calculated $(\overline{\alpha_1})_i$ one sets
 $(\overline{\alpha_2})_i=1-(\overline{\alpha_1})_i,$
\begin{equation}(\overline{E_1})_i=(\overline{F_1})_i-p_i^n(\overline{\alpha_1})_i,
(\overline{E_2})_i=(\overline{F_2})_i-p_i^n(\overline{\alpha_2})_i,\end{equation}
\begin{equation} (\overline{\rho_1})_i=\frac{(\overline{r_1})_i}{(\overline{\alpha_1})_i},(\overline{\rho_2})_i=\frac{(\overline{r_2})_i}{(\overline{\alpha_2})_i},\end{equation}
then the values of $\overline{p}_i$ obtained from (10,11) (by construction they are the same for $p_1$ and $p_2$). From  the formulas $\gamma_1=K_1-1,\gamma_2=K_2-1$ \cite{Munk,FFMunk} and from (8,9) one has the values $(\overline{v_\tau})_i$ and $(\overline{\Delta p})_i$  for $v_\tau$ and $\Delta p$ respectively.\\

$\bullet$ \textbf{fourth step: pressure correction.} Finally, we set
\begin{equation} (r_1v_1)_i^{n+1}=(\overline{r_1v_1})_i-\frac{r}{2}(\overline{\Delta p})_{i}((\overline{\alpha_1})_{i+1}-(\overline{\alpha_1})_{i-1})-\frac{r}{2}(\overline{\alpha_1})_{i}((\overline{p})_{i+1}-(\overline{p})_{i-1})+rhg(\overline{r_1})_i,\end{equation}
\begin{equation} (r_2v_2)_i^{n+1}=(\overline{r_2v_2})_i+\frac{r}{2}(\overline{\Delta p})_i((\overline{\alpha_1})_{i+1}-(\overline{\alpha_1})_{i-1})-\frac{r}{2}(\overline{\alpha_2})_i((\overline{p})_{i+1}-(\overline{p})_{i-1})+rhg(\overline{r_2})_i,\end{equation}
\begin{equation} (F_1^*)_i=(\overline{F}_1)_i-\frac{r}{2}(\overline{v_\tau})_i(\overline{\Delta p})_i((\overline{\alpha_1})_{i+1}-(\overline{\alpha_1})_{i-1})+rhg(\overline{r_1v_1})_i,\end{equation}
\begin{equation} (F_2^*)_i=(\overline{F}_2)_i+\frac{r}{2}(\overline{v_\tau})_{i}(\overline{\Delta p})_{i}((\overline{\alpha_1})_{i+1}-(\overline{\alpha_1})_{i-1})+rhg(\overline{r_2v_2})_{i}.\end{equation}
We set $(r_1)_i^{n+1}=(\overline{r_1})_i,\ (r_2)_i^{n+1}=(\overline{r_2})_i, (v_1)_i^{n+1}=\frac{(r_1v_1)_i^{n+1}}{(r_1)_i^{n+1}}$  and $(v_2)_i^{n+1}=\frac{(r_2v_2)_i^{n+1}}{(r_2)_i^{n+1}}$. \\ 

$\bullet$ \textbf{ fifth step: auxiliary algebraic calculations.} One has the values $(r_1)_i^{n+1}, (r_2)_i^{n+1}, (r_1v_1)_i^{n+1}, (r_2v_2)_i^{n+1}, (F_1)_i^*$ and $(F_2)_i^*$
as in situation (37), as well as the values $\overline{p}_i$ obtained in the third  step, that play here the role played by  $p_i^n$ in the third step. The formulas (39-42) with these values give  values of $(\alpha_j)_i^{n+1}$, then $ (E_j)_i^{n+1}$ and $(p_j)_i^{n+1}$ with $(p_1)_i^{n+1}=(p_2)_i^{n+1}$ from (42,43,10,11).
We denote this value $p_i^{n+1}=(p_1)_i^{n+1}=(p_2)_i^{n+1}$ and we set $(F_j)_i^{n+1}=(E_j)_i^{n+1}+p_i^{n+1}(\alpha_j)_i^{n+1}, j=1,2$.\\


\textbf{5. A numerical scheme for the system solved in time derivatives.}

A nicer formulation that involves time derivatives in the energy equations only in the variables $E_j, \ j=1,2$ has been obtained in \cite{Munk} p 2595-2596: the nonconservative temporal derivatives in (5,6) (formula (1) in \cite{Munk}) are eliminated (formulas (2,3) in \cite{Munk}). Equations (1-4) are unchanged. In \cite{FFMunk} p. 414   equations (5,6) are stated as

\begin{equation}\frac{\partial}{\partial t}(E_1)+\frac{\partial}{\partial x}(E_1v_1)
+(\alpha_1v_1-\eta\alpha_1\alpha_2(v_1-v_2))\frac{\partial p}{\partial x}+\eta\rho_2\alpha_1(c_2)^2\frac{\partial}{\partial x}(\alpha_1v_1+\alpha_2v_2))+v_{\tau}\tau_i=g\rho_1\alpha_1v_1,   \end{equation}

\begin{equation}\frac{\partial}{\partial t}(E_2)+\frac{\partial}{\partial x}(E_2v_2)
+(\alpha_2v_2+\eta\alpha_1\alpha_2(v_1-v_2))\frac{\partial p}{\partial x}+\eta\rho_1\alpha_2(c_1)^2\frac{\partial}{\partial x}(\alpha_1v_1+\alpha_2v_2))-v_{\tau}\tau_i=g\rho_2\alpha_2v_2,   \end{equation}
\\
where  
$c_j, \ j=1,2$ is the sound velocity
\begin{equation} (c_j)^2=\frac{p+K_jp_{\infty,j}}{\rho_j}\end{equation}
and where
\begin{equation} \eta=\frac{p}{\rho_1\alpha_2(c_1)^2+\rho_2\alpha_1(c_2)^2}.\end{equation}\\

It is proved in \cite{Munk} section 2.3.3 that the equations (48-51) are mathematically equivalent to (5-9) in the case of smooth solutions. Therefore it is important to check  this transformation  in the case of discontinuous solutions: a numerical verification will be done below in figures 1 and 2 in the hyperbolic and in the nonhyperbolic case respectively.\\


The numerical scheme is obtained as in section 2 from the three steps of transport, averaging and pressure correction, with auxiliary algebraic calculations. The transport step is made of equations 
(17-20) for $r_1, r_2, r_1v_1$  and $r_2v_2$ and the following  simple transport equations

\begin{equation} \frac{\partial E_j}{\partial t}+\frac{\partial (E_jv_j)}{\partial x}=0 \end{equation}
\\
for $E_1$ and $E_2$, that replace (21,22). Now we give  the scheme in more details. One starts with the values \begin{equation}\{(r_1)_i^n,(r_2)_i^n,(r_1v_1)_i^n,(r_2v_2)_i^n,(v_1)_i^n,(v_2)_i^n,(E_1)_i^n,(E_2)_i^n\}_{i\in\mathbb{Z}}.\end{equation} 

$\bullet$ \textbf{first step: transport.}
The values $\widetilde{r_1}, \widetilde{r_2}, \widetilde{r_1v_1}$ and $\widetilde{r_2v_2}$ are given by formulas (31-34). From (52) we obtain $\widetilde{E_j}, \ j=1,2$:
\begin{equation}(\widetilde{E_j})_i=(E_j)_i^n+r[(E_j)_{i-1}^n((v_j)_{i-1}^n)^+-(E_j)_{i}^n(|v_j|_{i}^n)+(E_j)_{i+1}^n((v_j)_{i+1}^n)^-].\end{equation}

$\bullet$ \textbf{second step: averaging.} The averaging step is identical to the one in section 2 for $r_1, r_2, r_1v_1$ and  $r_2v_2$; it is also done for $E_1$ and $E_2$. Then one states $(\overline{v_j})_i=\frac{(\overline{r_jv_j})_i}{(\overline{r_j})_i}, \ j=1,2$. At this point we have values $(\overline{r_1})_i, (\overline{r_2})_i, (\overline{r_1v_1})_i, (\overline{r_2v_2})_i, (\overline{v_1})_i, (\overline{v_2})_i, (\overline{E_1})_i$ and $(\overline{E_2})_i$. For the pressure correction term we need the values of $(\alpha_1)_i^n, (\alpha_2)_i^n, (\rho_1)_i^n, (\rho_2)_i^n$  and $p_i^n$. We obtain them  by auxiliary algebraic calculations in the third step.\\  

$\bullet$ \textbf{third step: auxiliary algebraic calculations.} Algebraic calculations with the values obtained in the second step, as this is done in section 2, would be more natural and would permit a better CFL condition ($r=0.002$ in the test of figure 1 instead of r=0.0012). Nevertheless we perform the algebraic calculations with the values (53) because this is more convenient for the proof of the theorem, and at the same time gives same numerical result. The value $\alpha_i^n=(\alpha_1)_i^n$ is obtained as the solution of the equation $AX^2+BX+C=0$ which lies in $]0,1[$, with 


\begin{equation}A=K_1p_{\infty,1}-K_2p_{\infty,2},\end{equation}
$ B=-(K_1-1)(E_1)_i^n-(K_2-1)(E_2)_i^n-K_1p_{\infty,1}+K_2p_{\infty,2}+\frac{1}{2}(K_1-1)(r_1v_1)_i^n(v_1)_i^n+$\begin{equation}\frac{1}{2}(K_2-1)(r_2v_2)_i^n(v_2)_i^n,\end{equation}
\begin{equation}C=(K_1-1)(E_1)_i^n-\frac{1}{2}(K_1-1)(r_1v_1)_i^n(v_1)_i^n.\end{equation}
\\
Then we set $(\alpha_2)_i^n=1-\alpha_i^n, (\rho_j)_i^n=\frac{(r_j)_i^n}{(\alpha_j)_i^n},$\\  
 \begin{equation}(p_j)_i^n=(K_j-1)(\frac{(E_j)_i^n}{(\alpha_j)_i^n}-\frac{1}{2}(\rho_j)_i^n((v_j)_i^n)^2)-K_j p_{\infty,j}, \  p_i^n=(p_1)_i^n=(p_2)_i^n,\end{equation} 
\begin{equation} 
 (\eta)_i^n=\frac{p_i^n}{(c_2)^2(\alpha_1)_i^n(\rho_2)_i^n+(c_1)^2(\alpha_2)_i^n(\rho_1)_i^n}.\end{equation}
We calculate the quantities $(v_\tau)_i^n$ and  $(\Delta p)_i^n$ from (8,9).\\


Then, dropping the superscript $n$ for convenience in second members of (60,61), we set (formulas 2.10,2.11 in \cite{FFMunk} p. 414)\\
$$ (T_1)_i^n=(\alpha_1)_i(v_1)_i(p_{i+1}-p_{i-1})-\eta_i(\alpha_1)_i(\alpha_2)_i((v_1)_{i}-(v_2)_{i})(p_{i+1}-p_{i-1})+$$  $$\eta_i(\rho_2)_i(\alpha_1)_i((c_2)^2)_i((\alpha_1)_{i+1}(v_1)_{i+1}-(\alpha_1)_{i-1}(v_1)_{i-1})+\eta_i(\rho_2)_i(\alpha_1)_i((c_2)^2)_i$$\begin{equation}((\alpha_2)_{i+1}(v_2)_{i+1}-(\alpha_2)_{i-1}(v_2)_{i-1})+(v_\tau)_i(\Delta p)_i((\alpha_1)_{i+1}-(\alpha_1)_{i-1}),\end{equation}

$$ (T_2)_i^n=(\alpha_2)_i(v_2)_i(p_{i+1}-p_{i-1})+\eta_i(\alpha_1)_i(\alpha_2)_i((v_1)_{i}-(v_2)_{i})(p_{i+1}-p_{i-1})+$$  $$\eta_i(\rho_1)_i(\alpha_2)_i((c_1)^2)_i((\alpha_1)_{i+1}(v_1)_{i+1}-(\alpha_1)_{i-1}(v_1)_{i-1})+\eta_i(\rho_1)_i(\alpha_2)_i((c_1)^2)_i$$\begin{equation}((\alpha_2)_{i+1}(v_2)_{i+1}-(\alpha_2)_{i-1}(v_2)_{i-1})-(v_\tau)_i(\Delta p)_i((\alpha_1)_{i+1}-(\alpha_1)_{i-1}).\end{equation}\\


$\bullet$ \textbf{fourth step: pressure correction.} For $r_1v_1$ and $r_2v_2$ one uses formulas (44) and (45) with the values of $(\Delta p)_i^n,\alpha_i^n,p_i^n$ obtained in the third step and the values $(r_j)_i^{n}$. For $E_1$ and $E_2$, we set
\begin{equation} (E_j)_i^{n+1}=(\overline{E_j})_i^n-\frac{1}{2}r (T_j)_i^n+rhg(r_jv_j)_i^n, \ j=1,2.\end{equation}
At this point we have values $(r_1)_i^{n+1}, (r_2)_i^{n+1}, (r_1v_1)_i^{n+1}, (r_2v_2)_i^{n+1}$, $(v_1)_i^{n+1}=\frac{(r_1v_1)_i^{n+1}}{(r_1)_i^{n+1}}, (v_2)_i^{n+1}=\frac{(r_2v_2)_i^{n+1}}{(r_2)_i^{n+1}}, (E_1)_i^{n+1}$ and $(E_2)_i^{n+1}$, which permits to continue the induction on $n$.\\ \\


\textbf{6. Numerical Results on the Toumi shock tube problem.}\\

INSERT FIGURE 1\\

\textit{Figure 1. Comparison of the hyperbolic models ($\delta=2$) solved (black) and nonsolved (red, nearly completely hidden under the black) in time derivative with the  scheme of sections 3 and 4:  exact coincidence.}\\

 The canonical nonconservative form (not solved in time derivative)  and the rewritten model (solved in time derivative) have same solution  in presence of the additional term that ensures hyperbolicity \cite{Munk} p. 2595 (i.e. $\delta=2$ for both models). Here $t= 0.06, 10.000$ space steps, $r=0.0012$ and  $a=0.3$ for both systems. Further, this solution is identical to the one presented in \cite{Munk} pp. 2615,2617 and in  \cite{FFMunk}  p. 437. \\

INSERT FIGURE 2\\

\textit{Figure 2. Comparison of the nonhyperbolic models ($\delta=0)$ solved (black) and nonsolved (red, nearly completely  hidden under the black) in time derivative with the scheme of sections 3 and 4:  exact coincidence.}\\

 The canonical nonconservative form (not solved in time derivative)  and the rewritten model  (solved in time derivative) have same solution  in absence of the additional term that ensures hyperbolicity  (i.e. $\delta=0$ for both models). Here $t=0.06, 4000$ space steps, $r=0.0009$ and  $a=0.3$ for both systems. \\

The results observed in figures 1 and 2 constitute a numerical verification of our scheme and of the validity of the transformation of equations in case of discontinuous solutions \cite {Munk}. Figure 2 also brings a numerical confirmation that some nonlinear nonhyperbolic systems could have well posed solutions as noticed by B. Keyfitz et al. in  
\cite{ Keyfitz0,Keyfitzrev1,Keyfitzrev2,Keyfitzrev3,Keyfitzlivro,Keyfitzart1,Keyfitzart2,Keyfitzart3,KeyfitzLopes,KeyfitzSan}.\\

INSERT FIGURE 3\\

\textit{Figure 3. Comparison of the hyperbolic model ($\delta=2$, black) and the nonhyperbolic model ($\delta=0$, red)   on 2000 space steps with the scheme of section 4.}\\ 

 Figure 3 shows numerical solutions of the system solved in time derivative with and without additional term at $t=0.06, 2000$ space steps, $r=0.0012$ and $a=0.3$ for both  systems. One observes that the results are really close, in particular the step values are always the same, with an irregularity in gas volume fraction and peaks in  gas velocity and liquid velocity in absence of the additional term (i.e. $\delta=0$). After a  longer time-with correspondingly a longer tube-we have observed that the aspects of both solutions did not change significantly.\\


The first observation on figures 1 and 2 is that the two forms of the system, solved and nonsolved in time derivative,  do give exactly the same numerical solutions even in presence of shocks and other irregularities. Then we observe that the hyperbolic form of the system gives   the solution obtained in \cite{Munk} pp. 2615,2617 and \cite{FFMunk}  p. 437 by  completely different numerical methods.
We observe on figure 3  differences due to the additional term that renders the system hyperbolic, \cite{Munk} p. 2595, but also there are strong similarities such as all step values and three physical variables practically identical: pressure, gas temperature and liquid temperature. What is the status of the irregularity observed in gas volume fraction, and of the peaks in gas velocity and liquid velocity in the nonhyperbolic model? They appear as a small  region with large gas and liquid velocities together with a liquid volume fraction close to one. First, one observes that these results are to a large extent  independent of the values of the space steps, of the number of iterations and of the averaging. The top values and width of the peaks are well defined; they  grow  slowly  with time. They are reproduced from  the higher order scheme  even with a small number of space steps as expected for a genuine solution (figure 5). Therefore they do not look like numerical artefacts.  The point to decide if these peaks, which could be easily observed by engineers, are really present  in the physical solution  depends on experimental measurements of the gas kick phenomenon.\\

 Experimental reproductions of the gas kick  and measurements have been done using a vertical 1240 meters deep test well \cite{Avelar}. First we note that when we reproduce all numerical tests  in this paper in a vertical (100 meters long) tube (instead of the horizontal 100 meters long tube under consideration) one observes that all results are practically unchanged, and also in presence of viscosity. Therefore, the results of the nonhyperbolic model in a vertical tube still show the same  peak in liquid velocity conjointly with a liquid volume fraction close to one (top left figure on the gas volume fraction), which permits some attempts of comparison with the experimental results. Experimental measurements of the liquid flow rate ($m^3/s$) are shown  in \cite{Avelar}, figure 4 p. 17.  Experimentalists observed  a peak in liquid flow rate far higher than the surrounding step values  (about 2 or 3 times higher: \cite{Avelar} figure 4 p. 17),  which was underestimated   by the hyperbolic model they used, \cite{Avelar}  p.15. Although these experimental conditions are somewhat different from the Toumi numerical test presented in this paper, even in a vertical tube and taking viscosity into account,  this experimental observation and the analogy with the results in figure 3 suggest that the peak in liquid flow rate observed  in the nonhyperbolic model  could perhaps be physically realistic and, perhaps,  be, in part, interpreted in   the gas kick phenomenon.\\

 Viscosity is present in the physical problem (\cite{Evje1,Evje2,Evje3,EvjeFlattenFriis,EvjeFriis1,EvjeFriis2,EvjeKarlsen1,EvjeKarlsen2}). An easy modification of the scheme by introduction of the viscosity into the transport step  with the usual discretization of second order derivatives    in the two momentum equations and the two energy equations with various relevant viscosity coefficients for gas (about $10^{-3}$) and liquid (from $10^{-3}$ to $10^2$) did not show a significant influence of the presence of viscosity  both for the hyperbolic and the nonhyperbolic versions. 
 Now we compare the  scheme of order 3 in space exposed in section 6 with the order 1 scheme used up to now.\\

INSERT FIGURE 4\\

\textit{Figure 4. Hyperbolic model ($\delta=2$). Comparison of the order 3 scheme   in space (black) and of  the order 1  scheme of section 4 on 2000 space steps.}\\

 We observe that the $p=3$ scheme (black color,  section 6) gives exactly the same result as the previous (i.e. $p=1$) scheme  (red color) but with a better accuracy since one observes details in figure 1 and $\cite{Munk}$ which do not appear from the  previous scheme  (red) with same value of the space step; 2000 space steps, $r=0.001$ and $a=0.2$ for the previous  scheme, $r=0.0002$ and $a=0.025$ for the $p=3$ scheme\\


INSERT FIGURE 5\\

\textit{Figure 5. Nonhyperbolic model ($\delta=0$). Comparison of the order 3 scheme in space (black) and of the order 1 scheme  of section 3 (red) on 500 space steps. The result from this last scheme is  far from the final results  in figures 2 and 3 while the scheme of section 7 gives  this final result.}\\

 We observe that with 500 space steps the $p= 3$ scheme (black) gives practically  the result of the $p=1$ scheme with 2000 cells (figure 3), with further some minor parasite oscillations; $r=0.001$ and $a=0.2$ for the $p=1$ scheme, $r=10^{-5}$ and $a=0.06$ for the $p=3$ scheme. With 500 space steps the $p=1$ scheme of section 3 does not give a correct result since the peak in gas velocity is not apparent. The observation that the peculiarities of the nonhyperbolic model in figures 2 and 3 appear  with a  better accuracy with the higher order $p=3$ scheme suggests again that they are really part of a solution of the nonhyperbolic model. \\


\textbf{7. Higher order  versions of the scheme.}  Since the system is demanding from the numerical viewpoint one faces the need of a more efficient scheme. In this section the scheme is easily transformed into higher order  schemes.  The scheme in use up to now is of order 1 in the time increment $\Delta t=rh$ and of order 1 in the space increment $\Delta x=h$. Since the value of $r$ is small the main point is to increase the order in space. The  scheme is made of 3 main steps: transport, averaging and pressure correction. We introduce a $p$-discretization in space , $p\in \mathbb{N},  \ p$  odd, which consists in performing all of these 3 steps using values $\omega_j,  i-p\leq j\leq i+p$ to calculate the next value $\omega_i$, which is very easy since the needed coefficients can be obviously calculated before the iterations by resolving two small $(p+1)\times(p+1)$ and $p\times p$ linear systems whose entries are powers of integers from Taylor's formula  so as to maximize the order in $h$ in the transport step and in the pressure correction  step. Explicit coefficients are given below for the  $p=3$ scheme; for $p=5$ the coefficients were numerically computed.  The scheme used up to now corresponds to $p=1$. Such a $p$-scheme will be of order  $  p$ in the space increment:  the transport step is discretized at order $p$ in $\Delta x$ and the other steps are discretized at order $>p$. Since the stability appeared bad already for $p=3$  (too small value of $r=\frac{\Delta t}{\Delta x}$), a second averaging step was used after the pressure correction step. Since this improvement was still unsatisfactory the two averaging steps and the transport step were performed implicitly. Then, in the case of ideal gases, the results became satisfactory with a value $r$ often comparable to the one  of the $p=1$ scheme in use up to now (an exception in figure 5 where $r$ is very small since the nonhyperbolic model is very demanding). In some cases when $p=3$ and $5$ the tests show parasite oscillations on a few meshes at the vertical of the discontinuities. They are essentially eliminated by a post-treatment on the final solution by  replacing the final values by mean values on a few cells on the left and on the right. The $p=3$ scheme and the similar  $p=5$ scheme were tested on the 1-D system of ideal gases and for the system in this paper. The $p=3$ scheme is used in figures 4 and 5 in this paper. Now we give the detailed formulas of the $p=3$ scheme. \\

$\bullet$  \textit{The transport step.} The values $\{\omega_i^n,v_i^n\}_{i\in\mathbb{Z}}$ are known and we compute  values improperly denoted $\{\omega_i^{n+1}\}_{i\in\mathbb{Z}}$ at the end of the transport step to be inserted into the next step of averaging. The implicit formula is:\\

$ \omega_i^{n+1}=\omega_i^{n}-r\{-\frac{1}{3}\omega_{i-3}^{n+1}v_{i-3}^{n,+}+\frac{3}{2}\omega_{i-2}^{n+1}v_{i-2}^{n,+}-3\omega_{i-1}^{n+1}v_{i-1}^{n,+}+\frac{11}{6}\omega_{i}^{n+1}v_{i}^{n,+}+$\begin{equation}\frac{11}{6}\omega_{i}^{n+1}v_{i}^{n,-}
-3\omega_{i+1}^{n+1}v_{i+1}^{n,-}+\frac{3}{2}\omega_{i+2}^{n+1}v_{i+2}^{n,-}-\frac{1}{3}\omega_{i+3}^{n+1}v_{i+3}^{n,-}\}.\end{equation}
Calculation is done with a sparse matrix having 7 nonzero diagonals. The sum of the $v^+$ terms (respectively of the $v^-$ terms) above replaces the  expression $\omega_{i-1}^{n}v_{i-1}^{n,+}-\omega_{i}^{n}v_{i}^{n,+}$ (resp. $-\omega_{i}^{n}v_{i}^{n,-}+\omega_{i+1}^{n}v_{i+1}^{n,-}$ ) of the explicit  $p=1$ scheme used up to now. Semi-implicit RK2 and RK4 methods in time can be used easily in place of the implicit Euler order one method (63), due to the independence in time and linearity in $\omega$ of the transport formula, and have given far better results.\\


$\bullet$  \textit{The two averaging  steps.} Starting from values improperly denoted $\{\omega_i^n\}_{i\in\mathbb{Z}}$ the implicit formula is (composition of two Laplacians):
\begin{equation} \omega_i^{n+1}=\omega_i^{n}+a\{-\omega_{i-3}^{n+1}+6\omega_{i-2}^{n+1}-15\omega_{i-1}^{n+1}+20\omega_{i}^{n+1}
-15\omega_{i+1}^{n+1}+6\omega_{i+2}^{n+1}-\omega_{i+3}^{n+1}\}.\end{equation}
The value $a$ is chosen from numerical tests as exposed for the order one scheme in \cite{ColombeauJCAM} p.17. One averaging step is done after the transport step and another averaging after the pressure correction step.\\

$\bullet$  \textit{The pressure correction   step.} We  use a centered discretization in space:  $\partial_x \omega$ is discretized by 

\begin{equation} \frac{1}{303}\omega_{i+3}-\frac{39}{404}\omega_{i+2}+\frac{69}{101}\omega_{i+1}-\frac{69}{101}\omega_{i-1}+\frac{39}{404}\omega_{i-2}-\frac{1}{303}\omega_{i-3}.\end{equation}\\
This can be done with the explicit Euler order one scheme in time (figures 4 and 5). It has been observed that the explicit RK2 and RK4 schemes in time (using the above discretization in space) give better results.\\

 The RK2 and RK4 discretizations in time  have not been presented in figures 4 and 5 (where we used only the Euler order one methods in time, implicit in transport and explicit in pressure correction) to stress the very significant improvements already obtained only from the space discretization which is original, since the various possible time discretizations are standard. \\

Extension to 2-D and 3-D can be done using the transport formulas in \cite{ColombeauNMPDE} pp. 96-100.
In 1-D we considered separately terms in $v^+$ and in $v^-$ in (63) that were discretized up to indices $k, \ i-p\leq k\leq i+p$. The $p=1$ expression involving $v^+$, respectively $v^-$, was transformed into the left, respectively right,  part of the $\{\dots\}$ in (63) for $p=3$. Similarly, in 2-D if the velocity vector is denoted by $(u,v)$ one considers separately terms in $(u^+,v^+), (u^+,v^-), (u^-,v^+)$  and  $ (u^-,v^-)$ in the $p=1$ transport formulas in \cite{ColombeauNMPDE} pp. 96-100, that we  replace for the calculation of $\omega_{i,j}^{n+1}$ in the $p$-scheme by expressions involving $\omega_{k,q}^n$ for $i-p\leq k\leq i+p, j-p\leq q\leq j+p$ chosen from Taylor's formula in 2-D  so as to maximize the order in the space step $\Delta x=\Delta y$. As an application one could study multifluid flows in multidimension as done in 1-D in this paper. The 2-D $p=1$ scheme was sufficient in \cite{ColombeauJCAM} to obtain exactly the results 
 presented in \cite{Lax1,Lax2} (obtained there by  order two schemes)  but in the more demanding context of multifluid flows  the higher order schemes could presumably be useful.\\


\textbf{8. Conclusion.} We have checked that a numerical scheme used in \cite{ColombeauSiam, ColombeauNMPDE, ColombeauJCAM, ColombeauJCAM2, ColombeauJDE} applies also for the standard system of one pressure model of two fluid flows with energy equations. The novelty is that our scheme works in the hyperbolic and nonhyperbolic cases, whether solved and nonsolved in time derivative,  that we obtain well defined approximate solutions that tend to satisfy the equations  somewhat independently on a choice of the nonconservative product   both in the hyperbolic case and in the nonhyperbolic case. Comparisons with previous works \cite{Munk, FFMunk} in the hyperbolic case show perfect agreement. In the nonhyperbolic case we observe the appearance of a peak in liquid flow rate, whose possible physical relevance in relation with the gas kick phenomenon remains to be elucidated. Our scheme has been explicitely stated at order three in space and possible extensions to higher order in multidimension are sketched.\\


\textit{Acknowledgements.}\\


\centering 
\includegraphics{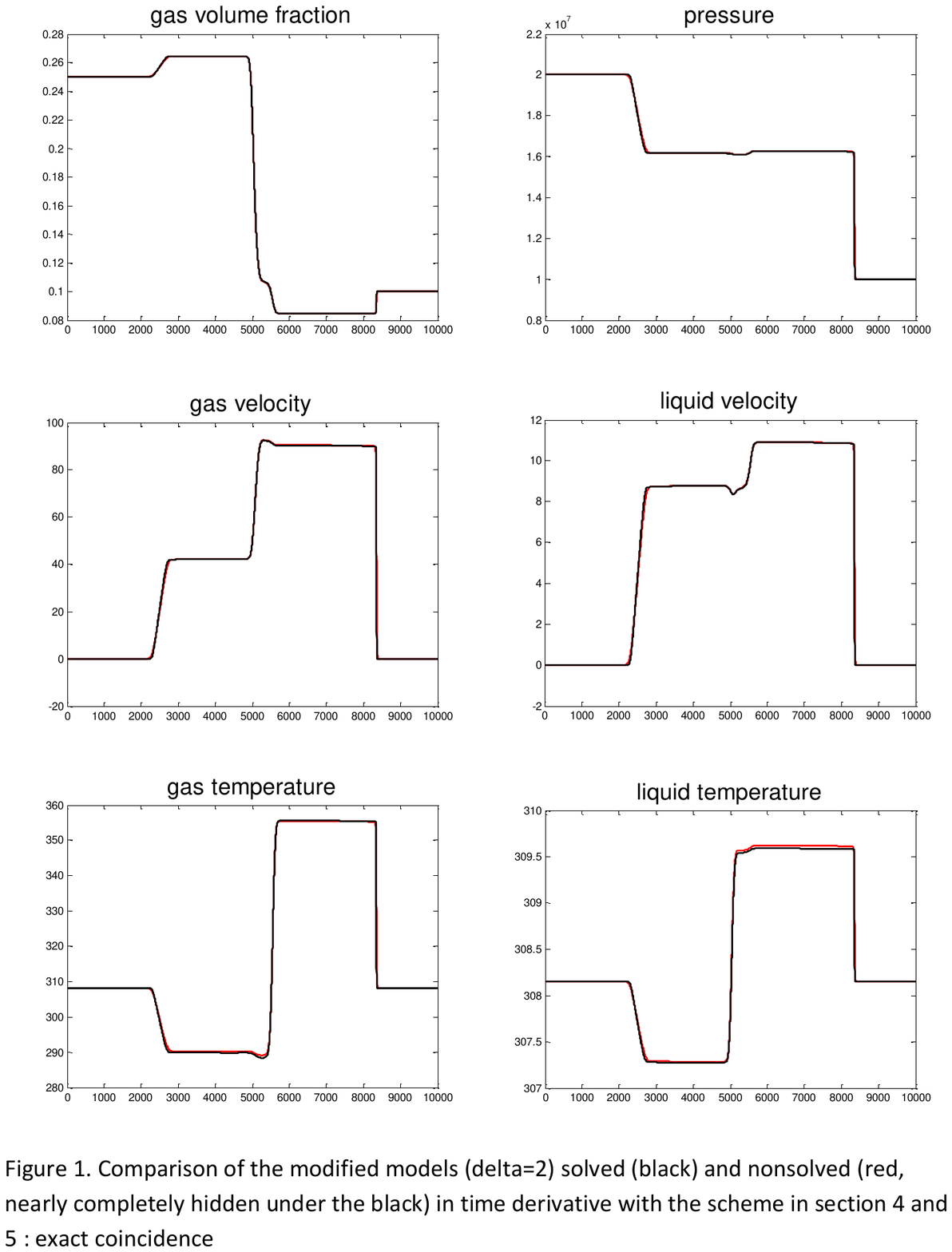}

\centering 
\includegraphics{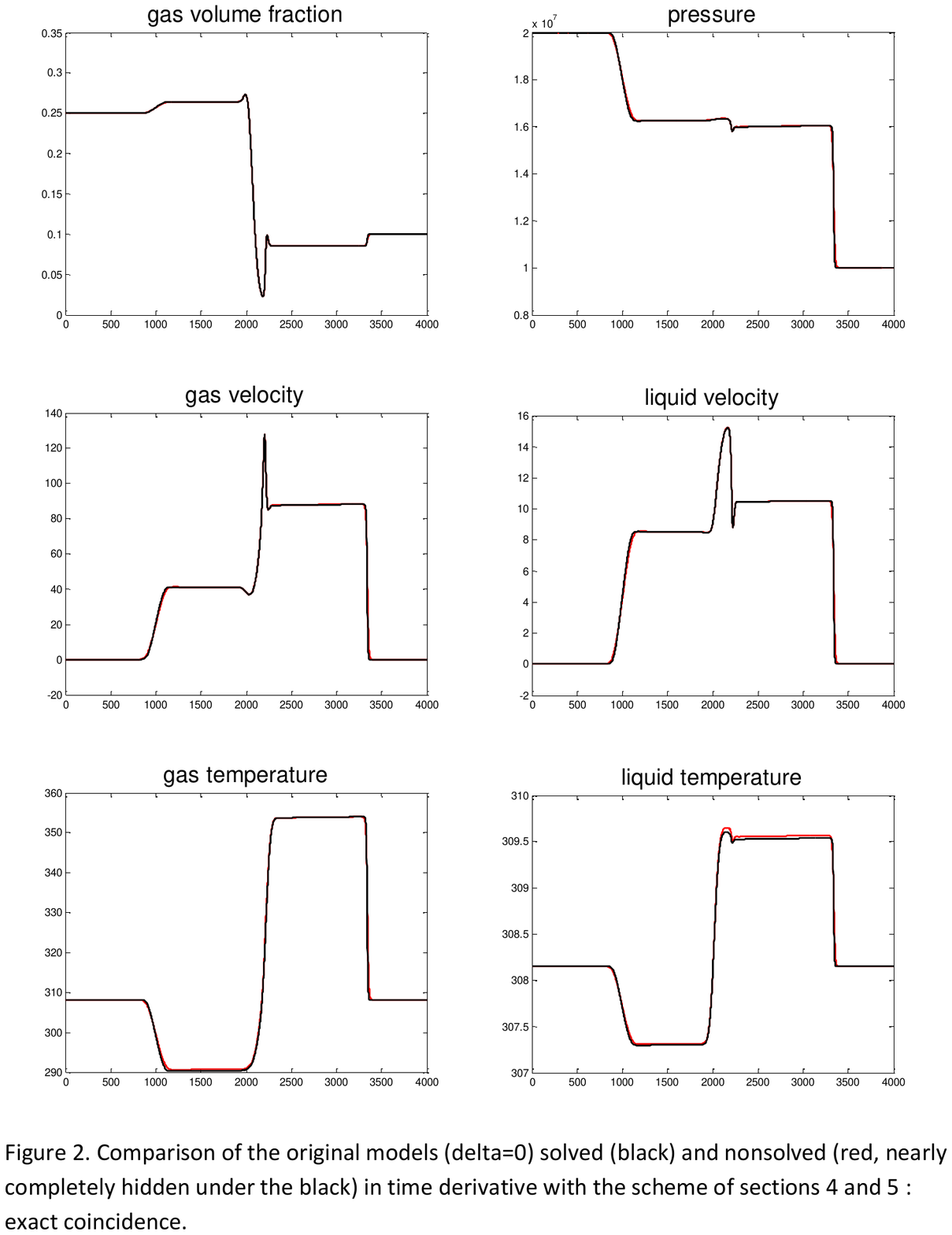}

\centering 
\includegraphics{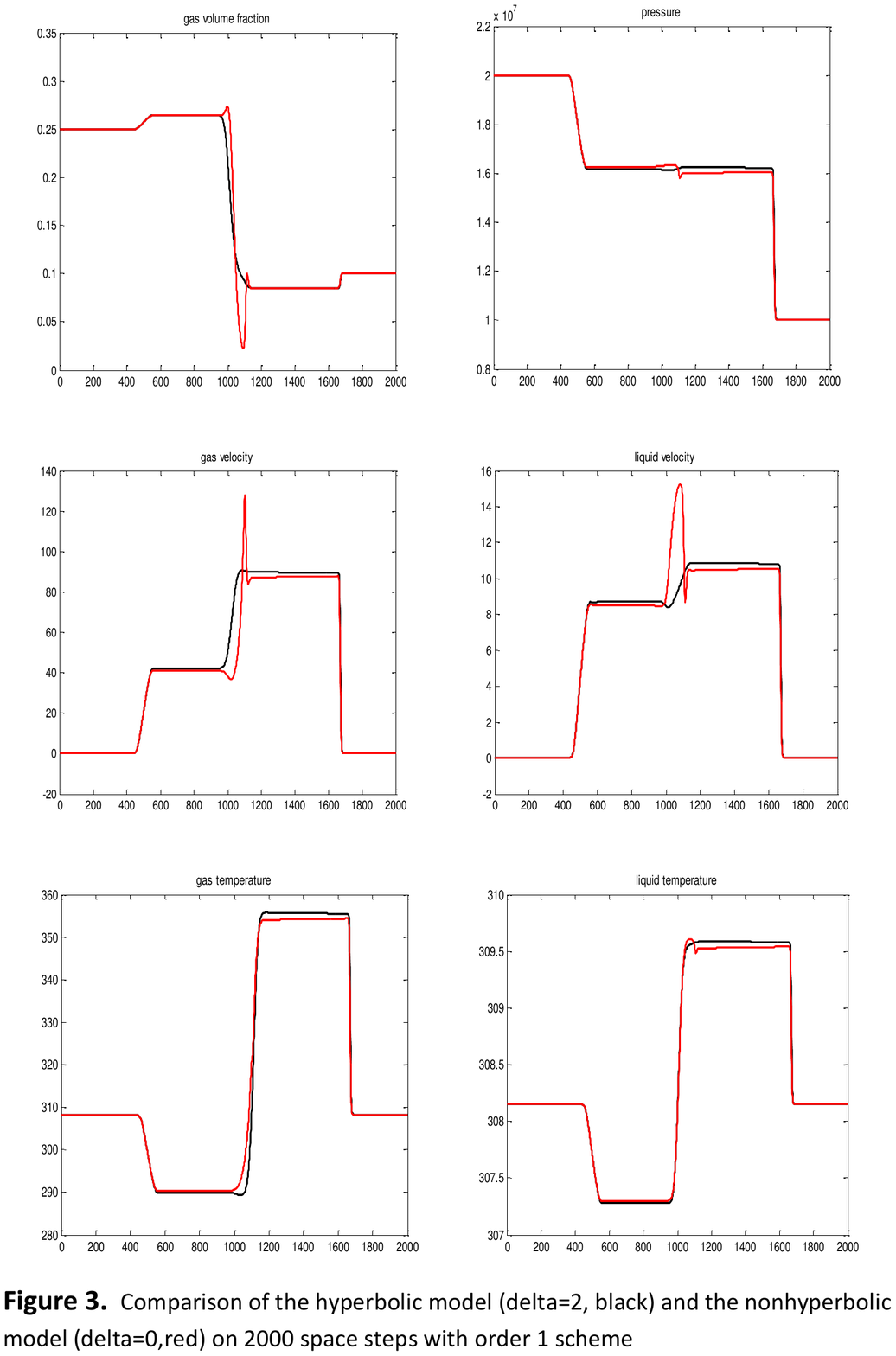}

\centering 
\includegraphics{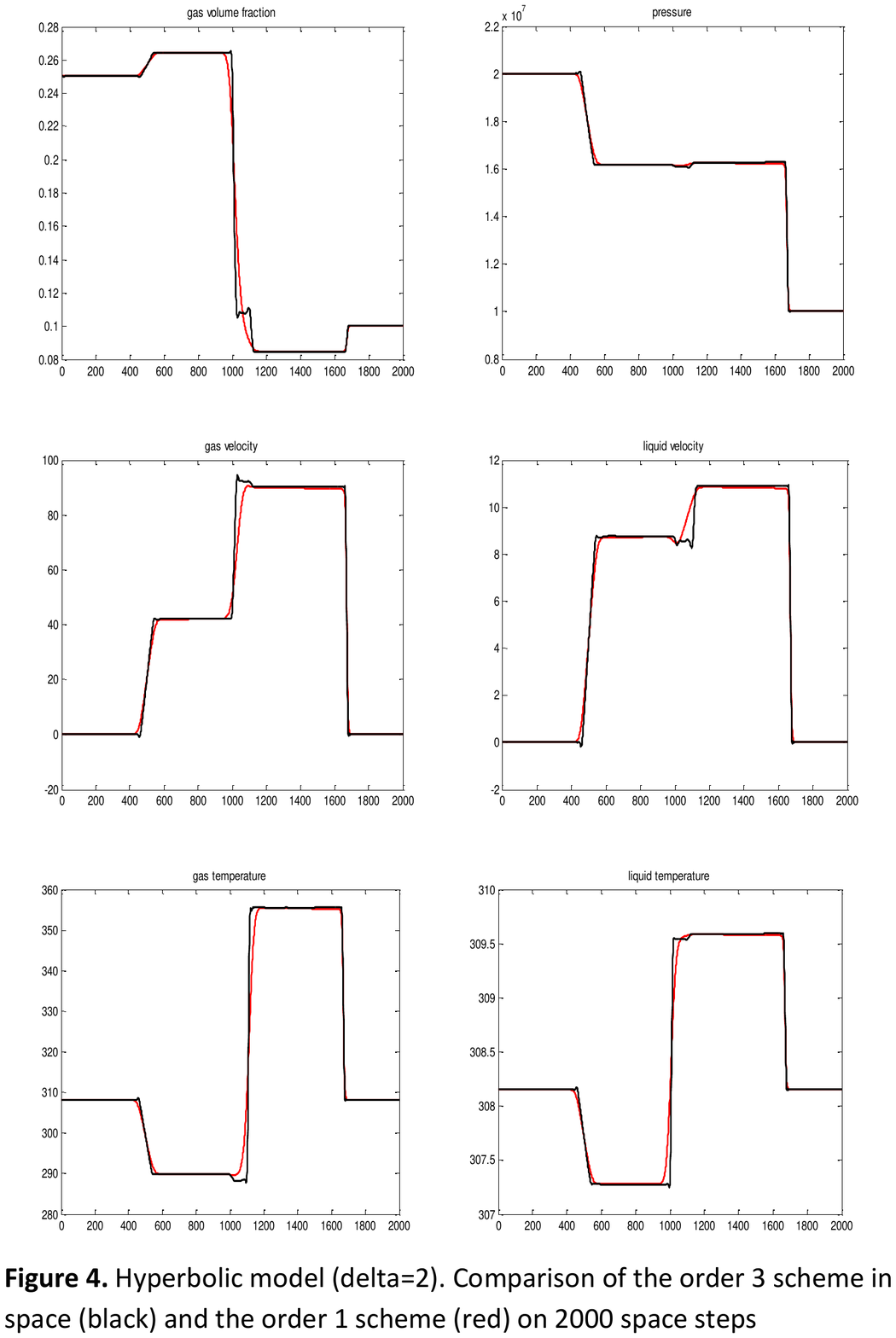}

\centering 
\includegraphics{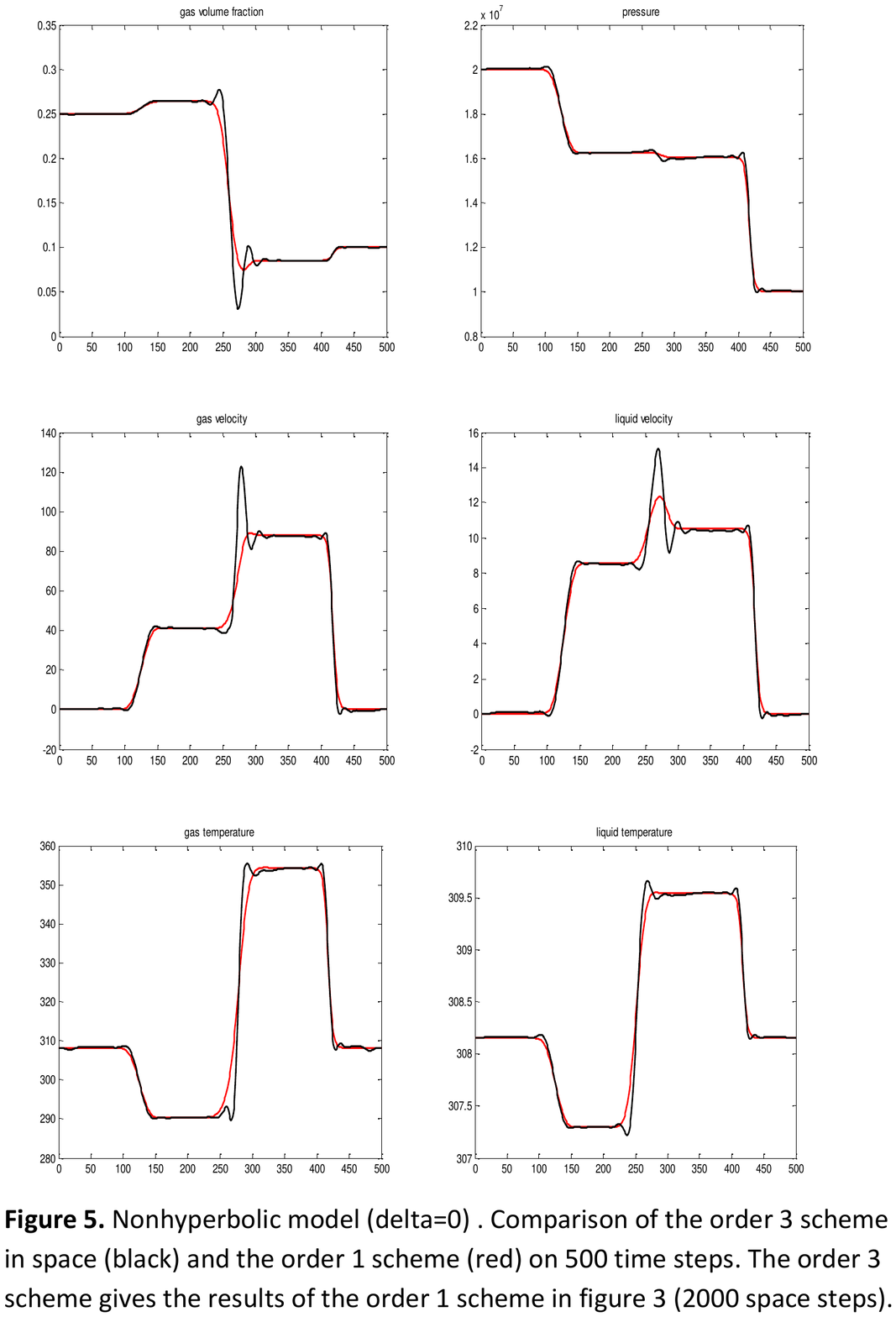}

 \end{document}